\documentclass{amsart}

\usepackage{amssymb}
\usepackage{amsmath}
\usepackage{amsthm}
\usepackage{amsfonts}

\usepackage{a4wide}
\usepackage{longtable}
\usepackage{multirow}

\usepackage{lineno}

\newtheorem{theorem}{Theorem}[section]

\newtheorem{corollary}{Corollary}[section]

\theoremstyle{definition}
\newtheorem{definition}{Definition}[section]

\newtheorem{example}{Example}[section]
\newtheorem{conjecture}{Conjecture}[section]

\begin{document}


\title{Normal lattice of certain metabelian \(p\)-groups \(G\) with \(G/G^\prime\simeq (p,p)\)}

\author{Daniel C. Mayer}
\address{Naglergasse 53\\8010 Graz\\Austria}
\email{algebraic.number.theory@algebra.at}
\urladdr{http://www.algebra.at}

\thanks{Research supported by the Austrian Science Fund (FWF): P 26008-N25}

\subjclass[2000]{Primary 20D15, 20F12, 20F14, secondary 11R37, 11R29, 11R11}
\keywords{lattice of normal subgroups, \(p\)-groups of derived length \(2\),
power-commutator presentations, central series, two-step centralizers,
second Hilbert \(p\)-class field, principalization of \(p\)-classes, quadratic fields}

\date{December 28, 2013}

\begin{abstract}
Let \(p\) be an odd prime.
The lattice of all normal subgroups and the terms of the lower and upper central series
are determined for all metabelian \(p\)-groups with generator rank \(d=2\)
having abelianization of type \((p,p)\) and minimal defect of commutativity \(k=0\).
It is shown that many of these groups are realized as Galois groups of second Hilbert \(p\)-class fields
of an extensive set of quadratic fields which are characterized by principalization types of \(p\)-classes.
\end{abstract}

\maketitle



\section{Introduction}
\label{s:Intro}

Let \(p\ge 3\) be an odd prime number,
and \(G=\langle x,y\rangle\) be a two-generated metabelian \(p\)-group
having an elementary bicyclic derived quotient \(G/G^\prime\) of type \((p,p)\).

Assume further that \(G\) is of order \(\lvert G\rvert=p^n\) with \(n\ge 2\),
and of nilpotency class \(\mathrm{cl}(G)=m-1\) with \(m\ge 2\).
Then \(G\) is of coclass \(\mathrm{cc}(G)=n-m+1=e-1\) with \(e\ge 2\).
Denote by
\[G=\gamma_1(G)>\gamma_2(G)=G^\prime>\ldots>\gamma_{m-1}(G)>\gamma_m(G)=1\]
the (descending) lower central series of \(G\), where \(\gamma_j(G)=\lbrack\gamma_{j-1}(G),G\rbrack\) for \(j\ge 2\), and by
\[1=\zeta_0(G)<\zeta_1(G)<\ldots<G^\prime=\zeta_{m-2}(G)<\zeta_{m-1}(G)=G\]
the (ascending) upper central series of \(G\), where \(\zeta_j(G)/\zeta_{j-1}(G)=\mathrm{Centre}(G/\zeta_{j-1}(G))\) for \(j\ge 1\).

Let \(s_2=t_2=\lbrack y,x\rbrack\) denote the main commutator of \(G\),
such that \(\gamma_2(G)=\langle s_2,\gamma_3(G)\rangle\).
By means of the two series
\(s_j=\lbrack s_{j-1},x\rbrack\) for \(j\ge 3\) and
\(t_\ell=\lbrack t_{\ell-1},y\rbrack\) for \(\ell\ge 3\)
of higher commutators and the subgroups
\(\Sigma_j=\langle s_j,\ldots,s_{m-1}\rangle\) with \(j\ge 3\) and \(T_\ell=\langle t_\ell,\ldots,t_{e+1}\rangle\) with \(\ell\ge 3\),
we obtain the following fundamental distinction of cases.

\begin{enumerate}
\item
The \textit{uniserial} case of a CF group (\textit{cyclic factors}) of coclass \(\mathrm{cc}(G)=1\) (maximal class),
where \(t_3\in\Sigma_3\), \(\gamma_3(G)=\langle s_3,\gamma_4(G)\rangle\), \(e=2\), and \(m=n\).
There are two subcases:\\
(1.1) \(t_3=1\in\gamma_m(G)\), where \(G\) contains an abelian maximal subgroup and \(k=0\),\\
(1.2) \(1\ne t_3\in\gamma_{m-k}(G)\), \(1\le k\le m-4\), where all maximal subgroups are non-abelian.
\item
The \textit{biserial} case of a non-CF or BCF group (\textit{bicyclic or cyclic factors}) of coclass \(\mathrm{cc}(G)\ge 2\),
where \(t_3\not\in\Sigma_3\), \(\gamma_3(G)=\langle s_3,t_3,\gamma_4(G)\rangle\), \(e\ge 3\), and \(m<n\).
Again there exist two subcases, characterized by the \textit{defect of commutativity} \(k\) of \(G\):\\
(2.1) \(t_{e+1}=1\in\gamma_m(G)\), where \(\Sigma_3\cap T_3=1\) and \(k=0\),\\
(2.2) \(1\ne t_{e+1}\in\gamma_{m-k}(G)\), for some \(k\ge 1\), where \(\Sigma_3\cap T_3\le\gamma_{m-k}(G)\).
\end{enumerate}

In this article, we are interested in two-generator metabelian \(p\)-groups \(G=\langle x,y\rangle\)
of coclass \(\mathrm{cc}(G)\ge 2\) having the convenient property \(\Sigma_3\cap T_3=1\), resp. \(k=0\),
where the product \(\Sigma_3\times T_3\) is direct and coincides with the major part of the \textit{normal lattice} of \(G\),
as shown in Figure
\ref{fig:NormalLatticeFigure}.

\begin{definition}
\label{dfn:Diamond}
A pair \((U,V)\) of normal subgroups of a \(p\)-group \(G\),
such that \(V<U\le G\) and \((U:V)=p^2\),
is called a \textit{diamond} if the quotient \(U/V\) is abelian of type \((p,p)\).
\end{definition}

If \((U,V)\) is a diamond and \(U=\langle u_1,u_2,V\rangle\),
then the \(p+1\) intermediate subgroups of \(G\) between \(U\) and \(V\) are given by
\(\langle u_2,V\rangle\) and \(\langle u_1u_2^{i-2},V\rangle\) with \(2\le i\le p+1\).

\newpage


\section{The normal lattice}
\label{s:NormalLattice}

In this section, let \(G=\langle x,y\rangle\) be a metabelian \(p\)-group
with two generators \(x,y\), having abelianization \(G/G^\prime\) of type \((p,p)\)
and satisfying the independence condition \(\Sigma_3\cap T_3=1\), that is,
\(G\) is a metabelian \(p\)-group with defect of commutativity \(k=0\)
\cite[\S\ 3.1.1, p. 412, and \S\ 3.3.2, p. 429]{Ma}.
We assume that \(G\) is of coclass \(\mathrm{cc}(G)\ge 2\),
since the normal lattice of \(p\)-groups of maximal class
has been determined by Blackburn \cite{Bl}.

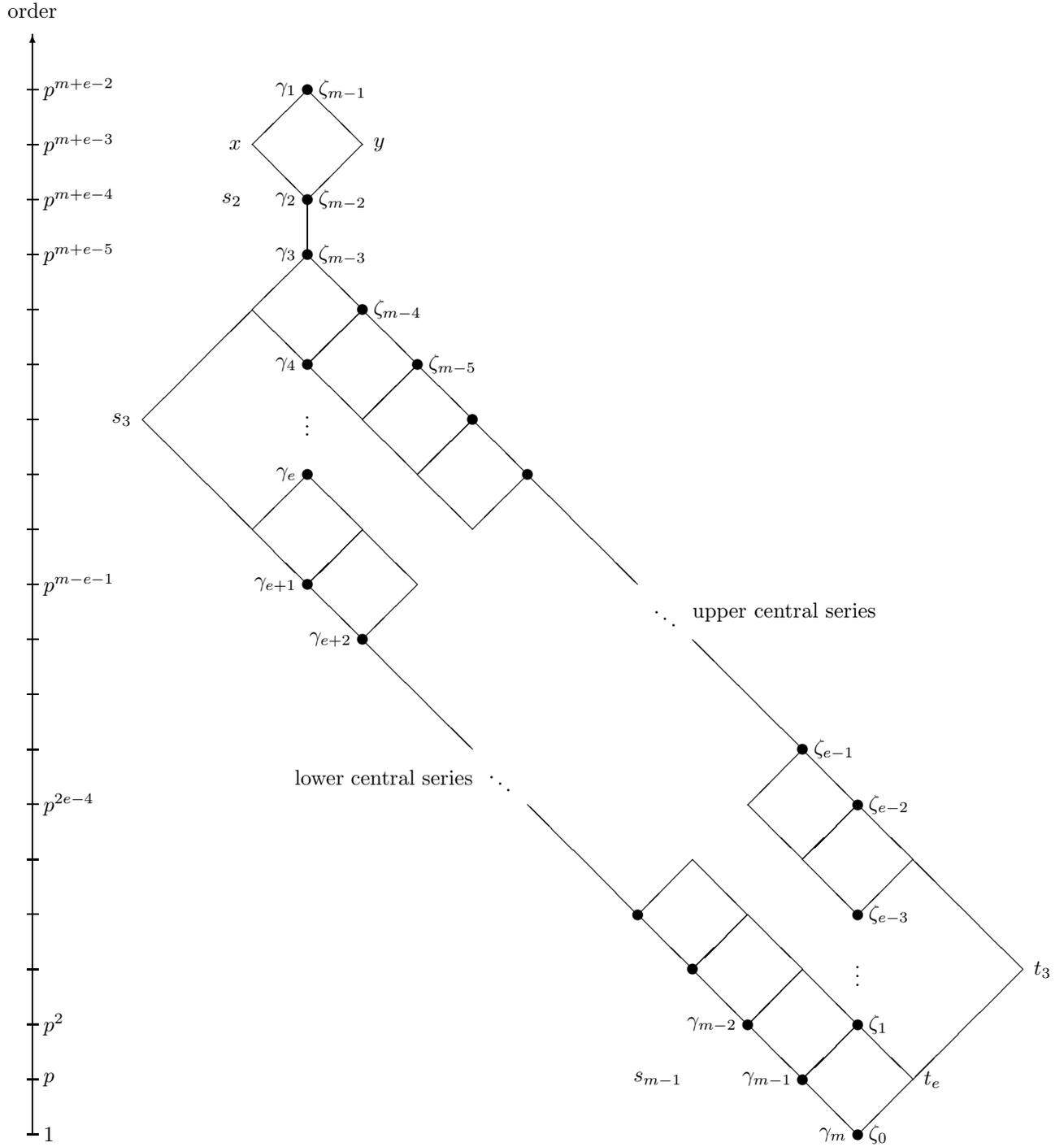
\begin{figure}[ht]
\caption{Full normal lattice, including lower and upper central series, of a \(p\)-group \(G\) with \(G/G^\prime\simeq (p,p)\), \(\mathrm{cl}(G)=m-1\), \(\mathrm{cc}(G)=e-1\), \(\mathrm{dl}(G)=2\), \(k(G)=0\).}
\label{fig:NormalLatticeFigure}


\setlength{\unitlength}{0.9cm}
\begin{picture}(18,22)(-5,-1)

\put(-5,20.3){\makebox(0,0)[cb]{order}}
\put(-5,19){\vector(0,1){1}}
\put(-4.8,19){\makebox(0,0)[lc]{\(p^{m+e-2}\)}}
\put(-4.8,18){\makebox(0,0)[lc]{\(p^{m+e-3}\)}}
\put(-4.8,17){\makebox(0,0)[lc]{\(p^{m+e-4}\)}}
\put(-4.8,16){\makebox(0,0)[lc]{\(p^{m+e-5}\)}}
\put(-4.8,10){\makebox(0,0)[lc]{\(p^{m-e-1}\)}}
\put(-4.8,6){\makebox(0,0)[lc]{\(p^{2e-4}\)}}
\put(-4.8,2){\makebox(0,0)[lc]{\(p^2\)}}
\put(-4.8,1){\makebox(0,0)[lc]{\(p\)}}
\put(-4.8,0){\makebox(0,0)[lc]{\(1\)}}
\multiput(-5.1,0)(0,1){20}{\line(1,0){0.2}}
\put(-5,0){\line(0,1){19}}


\put(-0.2,19){\makebox(0,0)[rc]{\(\gamma_1\)}}
\put(0,19){\circle*{0.2}}
\put(0.2,19){\makebox(0,0)[lc]{\(\zeta_{m-1}\)}}

\put(-1.2,18){\makebox(0,0)[rc]{\(x\)}}
\multiput(0,19)(-1,-1){2}{\line(1,-1){1}}
\multiput(0,19)(1,-1){2}{\line(-1,-1){1}}
\put(1.2,18){\makebox(0,0)[lc]{\(y\)}}

\put(-1.2,17){\makebox(0,0)[rc]{\(s_2\)}}
\put(-0.2,17){\makebox(0,0)[rc]{\(\gamma_2\)}}
\put(0,17){\circle*{0.2}}
\put(0.2,17){\makebox(0,0)[lc]{\(\zeta_{m-2}\)}}

\put(0,17){\line(0,-1){1}}

\put(-0.2,16){\makebox(0,0)[rc]{\(\gamma_3\)}}
\put(0,16){\circle*{0.2}}
\put(0.2,16){\makebox(0,0)[lc]{\(\zeta_{m-3}\)}}

\multiput(0,16)(-1,-1){2}{\line(1,-1){1}}
\multiput(0,16)(1,-1){2}{\line(-1,-1){1}}

\put(-0.2,14){\makebox(0,0)[rc]{\(\gamma_4\)}}
\put(0,14){\circle*{0.2}}

\put(1,15){\circle*{0.2}}
\put(1.2,15){\makebox(0,0)[lc]{\(\zeta_{m-4}\)}}

\multiput(1,15)(-1,-1){2}{\line(1,-1){1}}
\multiput(1,15)(1,-1){2}{\line(-1,-1){1}}

\put(2,14){\circle*{0.2}}
\put(2.2,14){\makebox(0,0)[lc]{\(\zeta_{m-5}\)}}

\multiput(2,14)(-1,-1){2}{\line(1,-1){1}}
\multiput(2,14)(1,-1){2}{\line(-1,-1){1}}

\put(3,13){\circle*{0.2}}

\multiput(3,13)(-1,-1){2}{\line(1,-1){1}}
\multiput(3,13)(1,-1){2}{\line(-1,-1){1}}

\put(4,12){\circle*{0.2}}

\put(4,12){\line(1,-1){2}}

\put(-3.2,13){\makebox(0,0)[rc]{\(s_3\)}}
\put(-3,13){\line(1,1){2}}
\put(-3,13){\line(1,-1){2}}

\put(0,13){\makebox(0,0)[cc]{\(\vdots\)}}

\put(-0.2,12){\makebox(0,0)[rc]{\(\gamma_e\)}}
\put(0,12){\circle*{0.2}}

\multiput(0,12)(-1,-1){2}{\line(1,-1){1}}
\multiput(0,12)(1,-1){2}{\line(-1,-1){1}}

\put(-0.2,10){\makebox(0,0)[rc]{\(\gamma_{e+1}\)}}
\put(0,10){\circle*{0.2}}

\multiput(1,11)(-1,-1){2}{\line(1,-1){1}}
\multiput(1,11)(1,-1){2}{\line(-1,-1){1}}

\put(0.8,9){\makebox(0,0)[rc]{\(\gamma_{e+2}\)}}
\put(1,9){\circle*{0.2}}

\put(1,9){\line(1,-1){2}}
\put(3.5,6.5){\makebox(0,0)[cc]{\(\ddots\)}}
\put(3,6.5){\makebox(0,0)[rc]{lower central series}}

\put(7,9.5){\makebox(0,0)[lc]{upper central series}}
\put(6.5,9.5){\makebox(0,0)[cc]{\(\ddots\)}}
\put(9,7){\line(-1,1){2}}

\put(9,7){\circle*{0.2}}
\put(9.2,7){\makebox(0,0)[lc]{\(\zeta_{e-1}\)}}

\multiput(9,7)(-1,-1){2}{\line(1,-1){1}}
\multiput(9,7)(1,-1){2}{\line(-1,-1){1}}

\put(10,6){\circle*{0.2}}
\put(10.2,6){\makebox(0,0)[lc]{\(\zeta_{e-2}\)}}

\multiput(10,6)(-1,-1){2}{\line(1,-1){1}}
\multiput(10,6)(1,-1){2}{\line(-1,-1){1}}

\put(10,4){\circle*{0.2}}
\put(10.2,4){\makebox(0,0)[lc]{\(\zeta_{e-3}\)}}

\put(13.2,3){\makebox(0,0)[lc]{\(t_3\)}}
\put(13,3){\line(-1,1){2}}
\put(13,3){\line(-1,-1){2}}

\put(10,3){\makebox(0,0)[cc]{\(\vdots\)}}

\put(6,4){\line(-1,1){2}}

\put(6,4){\circle*{0.2}}

\multiput(7,5)(-1,-1){2}{\line(1,-1){1}}
\multiput(7,5)(1,-1){2}{\line(-1,-1){1}}

\put(7,3){\circle*{0.2}}

\multiput(8,4)(-1,-1){2}{\line(1,-1){1}}
\multiput(8,4)(1,-1){2}{\line(-1,-1){1}}

\put(7.8,2){\makebox(0,0)[rc]{\(\gamma_{m-2}\)}}
\put(8,2){\circle*{0.2}}

\multiput(9,3)(-1,-1){2}{\line(1,-1){1}}
\multiput(9,3)(1,-1){2}{\line(-1,-1){1}}

\put(8.8,1){\makebox(0,0)[rc]{\(\gamma_{m-1}\)}}
\put(9,1){\circle*{0.2}}

\put(10,2){\circle*{0.2}}
\put(10.2,2){\makebox(0,0)[lc]{\(\zeta_1\)}}

\put(6.8,1){\makebox(0,0)[rc]{\(s_{m-1}\)}}
\multiput(10,2)(-1,-1){2}{\line(1,-1){1}}
\multiput(10,2)(1,-1){2}{\line(-1,-1){1}}
\put(11.2,1){\makebox(0,0)[lc]{\(t_e\)}}

\put(9.8,0){\makebox(0,0)[rc]{\(\gamma_m\)}}
\put(10,0){\circle*{0.2}}
\put(10.2,0){\makebox(0,0)[lc]{\(\zeta_0\)}}

\end{picture}

\end{figure}

\begin{theorem}
\label{thm:NormalLattice}
The complete normal lattice of \(G\)
contains
the heading diamond \((G,G^\prime)\)
and the rectangle \(\bigl((P_{j,\ell},P_{j+1,\ell+1})\bigr)_{3\le j\le m-1,\ 3\le\ell\le e}\) of trailing diamonds,
where \(P_{j,\ell}=\Sigma_j\times T_\ell\) for \(3\le j\le m\) and \(3\le\ell\le e+1\).
The structure of the normal lattice is visualized in Figure
\ref{fig:NormalLatticeFigure}.
\end{theorem}

Note that
\(P_{j,\ell}=\langle s_j,\ldots,s_{m-1}\rangle\times\langle t_\ell,\ldots,t_e\rangle=\langle s_j,t_\ell,P_{j+1,\ell+1}\rangle\)
for \(3\le j\le m-1\), \(3\le\ell\le e\).

\begin{conjecture}
\label{cnj:NormalLattice}
The complete normal lattice of \(G\)
consists exactly of the normal subgroups given in Theorem
\ref{thm:NormalLattice}.
\end{conjecture}

\begin{corollary}
\label{cor:NormalLattice}
The total number of normal subgroups of \(G\) is given by
\[me-(m+2e)+6+\lbrack me-(2m+3e)+7\rbrack\cdot (p-1),\]
in particular, for \(p=3\) it is given by
\[3me-(5m+8e)+20.\]
\end{corollary}

\begin{corollary}
\label{cor:CentralSeries}
Blackburn's two-step centralizers of \(G\)
\cite{Bl}
are given by

\begin{equation*}
\chi_j(G)=
\begin{cases}
G^\prime                  \text{ for } 1\le j\le e-1,\\
\langle y,G^\prime\rangle \text{ for } e\le j\le m-2,\\
G                         \text{ for } j\ge m-1,
\end{cases}
\end{equation*}

\noindent
in particular, none of the maximal subgroups of \(G\)
occurs as a two-step centralizer, when \(e=m-1\).

\begin{enumerate}

\item
The factors of the lower central series of \(G\) are given by

\begin{equation*}
\gamma_j(G)/\gamma_{j+1}(G)\simeq
\begin{cases}
(p,p) \text{ for } j=1 \text{ and } 3\le j\le e,\\
(p)   \text{ for } j=2 \text{ and } e+1\le j\le m-1.
\end{cases}
\end{equation*}

\item
The terms of the lower central series of \(G\) are given by

\begin{equation*}
\gamma_j(G)=
\begin{cases}
\langle x,y,G^\prime\rangle    \text{ for } j=1,\\
\langle s_2,\gamma_3(G)\rangle \text{ for } j=2,\\
P_{j,j}                        \text{ for } 3\le j\le e,\\
\Sigma_j                       \text{ for } e+1\le j\le m-1.
\end{cases}
\end{equation*}

\item
The factors of the upper central series of \(G\) are given by

\begin{equation*}
\zeta_j(G)/\zeta_{j-1}(G)\simeq
\begin{cases}
(p,p) \text{ for } 1\le j\le e-2 \text{ and } j=m-1,\\
(p)   \text{ for } e-1\le j\le m-2.
\end{cases}
\end{equation*}

\item
The terms of the upper central series of \(G\) are given by

\begin{equation*}
\zeta_j(G)=
\begin{cases}
P_{m-j,e+1-j}                     \text{ for } 1\le j\le e-2,\\
P_{m-j,3}                         \text{ for } e-1\le j\le m-3,\\
\langle s_2,\zeta_{m-3}(G)\rangle \text{ for } j=m-2,\\
\langle x,y,\zeta_{m-2}(G)\rangle \text{ for } j=m-1.
\end{cases}
\end{equation*}

\end{enumerate}

\end{corollary}

\newpage


\begin{proof}

We prove the invariance of all claimed normal subgroups
under inner automorphisms of \(G=\langle x,y\rangle\).

It is well known that the subgroups in the heading diamond are normal,
since they contain the commutator subgroup \(G^\prime=\gamma_2(G)\).

We start the proof with the tops of trailing diamonds.
For \(g\in P_{j,\ell}\) and \(s\in G^\prime\) we have
\(s^{-1}gs=s^{-1}sg=g\),
since \(P_{j,\ell}<G^\prime\), for \(j\ge 3\), \(\ell\ge 3\),
and \(G\) was assumed to be metabelian.
Now, \(P_{j,\ell}\) is the direct product of \(\Sigma_j\) and \(T_\ell\),
since we suppose that \(\Sigma_3\cap T_3=1\).
So it suffices to show invariance of \(\Sigma_j\) and \(T_\ell\)
under conjugation with the generators \(x\) and \(y\) of \(G\).
We have
\(x^{-1}s_jx=s_j\lbrack s_j,x\rbrack=s_js_{j+1}\in\Sigma_j\) and
\(y^{-1}s_jy=s_j\lbrack s_j,y\rbrack=s_j\in\Sigma_j\) for \(j\ge 3\).
And similarly we have
\(x^{-1}t_\ell x=t_\ell\lbrack t_\ell,x\rbrack=t_\ell\in T_\ell\) and
\(y^{-1}t_\ell y=t_\ell\lbrack t_\ell,y\rbrack=t_\ell t_{\ell+1}\in T_\ell\)
for \(\ell\ge 3\).

Next we prove invariance of intermediate groups between top and bottom of trailing diamonds.
They are of the shape
\(\langle t_\ell,P_{j+1,\ell+1}\rangle\) or
\(\langle s_jt_\ell^i,P_{j+1,\ell+1}\rangle\) with \(0\le i\le p-1\).
For \(t_\ell\), invariance has been shown above. So we investigate \(s_jt_\ell^i\).
We have \(x^{-1}s_jt_\ell^ix=x^{-1}s_jx(x^{-1}t_\ell x)^i=s_js_{j+1}t_\ell^i\),
where \(s_{j+1}\in P_{j+1,\ell+1}\), and
\(y^{-1}s_jt_\ell^iy=y^{-1}s_jy(y^{-1}t_\ell y)^i=s_jt_\ell^it_{\ell+1}^i\),
where \(t_{\ell+1}^i\in P_{j+1,\ell+1}\).
(Here we probably are tacitly using power conditions like
\(s_j^p\in\Sigma_{j+1}\) for \(j\ge 3\) and
\(t_\ell^p\in T_{\ell+1}\) for \(\ell\ge 3\).)

Thus we have proved the invariance of all claimed normal subgroups
under inner automorphisms.

\bigskip
The number of all (heading and trailing) diamonds of the normal lattice is
\(1+(m-1-2)\cdot (e-2)
=1+(m-3)\cdot (e-2)
=1+me-2m-3e+6
=me-(2m+3e)+7\).

There are \(p-1\) inner vertices of valence \(2\) in each diamond,
which gives a total of\\
\((me-\lbrack 2m+3e\rbrack+7)\cdot (p-1)\) inner vertices.

The remaining (outer) vertices form the heading square and the trailing rectangle with\\
\(4+(m-1+1-2)\cdot (e+1-2)
=4+(m-2)\cdot (e-1)
=4+me-m-2e+2
=me-(m+2e)+6\)
vertices.

Outer and inner vertices together form a lattice of
\(me-(m+2e)+6+(me-\lbrack 2m+3e\rbrack+7)\cdot (p-1)\)
normal subgroups.

For \(p=3\), this formula yields
\(me-m-2e+6+2me-4m-6e+14
=3me-(5m+8e)+20\).

\bigskip
For each \(j\ge 2\), Blackburn's two-step centralizer \(\chi_j(G)\) is defined as
the biggest intermediate group between \(G\) and \(G^\prime=\gamma_2(G)\) such that
\(\lbrack\gamma_j(G),\chi_j(G)\rbrack\le\gamma_{j+2}(G)\).
Since \(\lbrack\gamma_j(G),\gamma_2(G)\rbrack\le\gamma_{j+2}(G)\), for any \(j\ge 2\),
\(\chi_j(G)\) certainly contains \(\gamma_2(G)\).
Since
\(\lbrack s_j,x\rbrack=s_{j+1}\notin\gamma_{j+2}(G)\) for \(2\le j\le m-2\),
\(\lbrack t_\ell,y\rbrack=t_{\ell+1}\notin\gamma_{\ell+2}(G)\) for \(2\le\ell\le e-1\),
and \(e\le m-1\),
neither \(x\) nor \(y\) can be an element of \(\chi_j(G)\) for \(2\le j\le e-1\).
However, since
\(\lbrack t_e,y\rbrack=t_{e+1}=1\in\gamma_{e+2}(G)\) and
\(\lbrack s_e,y\rbrack=1\in\gamma_{e+2}(G)\),
we have \(\chi_j(G)=\langle y,\gamma_2(G)\rangle\) for \(e\le j\le m-2\),
provided that \(e\le m-2\).
Finally, since
\(\lbrack s_{m-1},x\rbrack=s_m=1\in\gamma_{m}(G)=\gamma_{m+1}(G)=1\),
the two-step centralizers \(\chi_j(G)\) with \(j\ge m-1\)
coincide with the entire group \(G\). 

\bigskip
The members of the lower central series can be constructed recursively by
\(\gamma_j(G)=\lbrack\gamma_{j-1}(G),G\rbrack\).
There is a unique ramification generating the series \(\Sigma_3\) and \(T_3\) for \(j=3\), since
\(\gamma_3(G)=\lbrack\gamma_2(G),G\rbrack=\lbrack\langle s_2,\gamma_3(G)\rangle,G\rbrack
=\langle\lbrack s_2,x\rbrack,\lbrack s_2,y\rbrack,\gamma_4(G)\rangle=\langle s_3,t_3,\gamma_4(G)\rangle\).
Otherwise the series \(\Sigma_3\) and \(T_3\) do not mix and we have
\(\gamma_j(G)=\lbrack\gamma_{j-1}(G),G\rbrack=\lbrack\langle s_{j-1},t_{j-1},\gamma_j(G)\rangle,G\rbrack\)\\
\(=\langle\lbrack s_{j-1},x\rbrack,\lbrack s_{j-1},y\rbrack,\lbrack t_{j-1},x\rbrack,\lbrack t_{j-1},y\rbrack,\gamma_{j+1}(G)\rangle
=\langle s_j,t_j,\gamma_{j+1}(G)\rangle\),
since \(\lbrack s_{j-1},y\rbrack=\lbrack t_{j-1},x\rbrack=1\) for \(j\ge 4\).
For \(j=e+1\) the bicyclic factors stop, since \(t_{e+1}=\lbrack t_e,y\rbrack=1\),
and \(\gamma_{e+1}\) is simply given by \(\Sigma_{e+1}\).

\bigskip
The members of the upper central series can be constructed recursively by
\(\zeta_j(G)/\zeta_{j-1}(G)=\mathrm{Centre}(G/\zeta_{j-1}(G))\).
All groups \(G\) with the assigned properties have a bicyclic centre
\(\zeta_1(G)=\langle s_{m-1},t_e\rangle\),
since
\(\lbrack s_{m-1},x\rbrack=\lbrack t_e,y\rbrack=1\).

Generally, the equations
\(\lbrack s_{m-j},x\rbrack=s_{m-(j-1)},\ \lbrack s_{m-j},y\rbrack=1,\ \lbrack t_{e+1-j},x\rbrack=1,\ \lbrack t_{e+1-j},y\rbrack=t_{e+1-(j-1)}\),
whose right sides are elements of \(\zeta_{j-1}(G)\),
show that \(s_{m-j}\) and \(t_{e+1-j}\) commute with all elements of \(G\) modulo \(\zeta_{j-1}(G)\).
Therefore, we have \(\zeta_j(G)=P_{m-j,e+1-j}\).

However, for \(j=e-1\) the bicyclic factors stop,
since \(\lbrack t_{e+1-j},x\rbrack=\lbrack t_2,x\rbrack=\lbrack s_2,x\rbrack=s_3\),
which is not contained in \(\zeta_{e-2}(G)\), except for \(e=m-1\).
Consequently, \(\zeta_j(G)=P_{m-j,3}\) for \(j\ge e-1\),
since it cannot contain \(t_2=s_2\).

\end{proof}

\newpage


\newpage


\section{Applications in Algebraic Number Theory}
\label{s:NumberTheoryApps}

Let \(K=\mathbb{Q}(\sqrt{D})\) be a quadratic number field with discriminant \(D\)
and denote by \(G=\mathrm{Gal}(\mathrm{F}_p^2(K)\vert K)\) the Galois group
of the second Hilbert \(p\)-class field \(\mathrm{F}_p^2(K)\) of \(K\),
that is, the maximal metabelian unramified \(p\)-extension of \(K\).
We recall that coclass and class of \(G\) are given by the equations
\(\mathrm{cc}(G)=r=e-1\) and \(\mathrm{cl}(G)=m-1\)
in terms of the invariants \(e\) and \(m\).
Due to our extensive computations for the papers
\cite{Ma0,Ma},
we are able to underpin the present theory of normal lattices
by numerical data concerning the \(2\,020\) complex and the \(2\,576\) real
quadratic fields with \(3\)-class group of type \((3,3)\)
and discriminant in the range \(-10^6<D<10^7\).

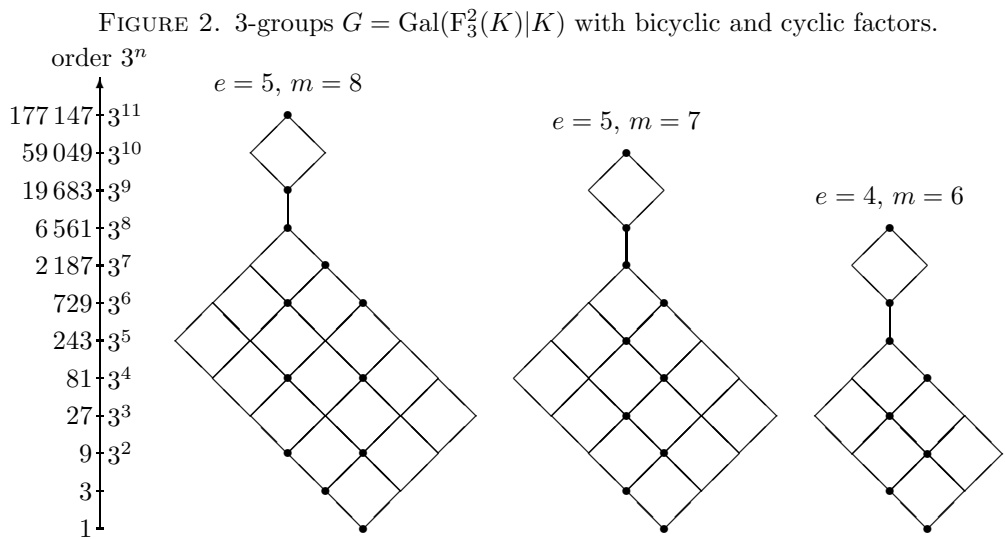
\begin{figure}[ht]
\caption{\(3\)-groups \(G=\mathrm{Gal}(\mathrm{F}_3^2(K)\vert K)\) with bicyclic and cyclic factors.}
\label{fig:BCFGroups}


\setlength{\unitlength}{0.5cm}
\begin{picture}(26,13)(-7,0)

\put(-5,12.3){\makebox(0,0)[cb]{order \(3^n\)}}
\put(-5,11){\vector(0,1){1}}
\put(-5.2,11){\makebox(0,0)[rc]{\(177\,147\)}}
\put(-4.8,11){\makebox(0,0)[lc]{\(3^{11}\)}}
\put(-5.2,10){\makebox(0,0)[rc]{\(59\,049\)}}
\put(-4.8,10){\makebox(0,0)[lc]{\(3^{10}\)}}
\put(-5.2,9){\makebox(0,0)[rc]{\(19\,683\)}}
\put(-4.8,9){\makebox(0,0)[lc]{\(3^9\)}}
\put(-5.2,8){\makebox(0,0)[rc]{\(6\,561\)}}
\put(-4.8,8){\makebox(0,0)[lc]{\(3^8\)}}
\put(-5.2,7){\makebox(0,0)[rc]{\(2\,187\)}}
\put(-4.8,7){\makebox(0,0)[lc]{\(3^7\)}}
\put(-5.2,6){\makebox(0,0)[rc]{\(729\)}}
\put(-4.8,6){\makebox(0,0)[lc]{\(3^6\)}}
\put(-5.2,5){\makebox(0,0)[rc]{\(243\)}}
\put(-4.8,5){\makebox(0,0)[lc]{\(3^5\)}}
\put(-5.2,4){\makebox(0,0)[rc]{\(81\)}}
\put(-4.8,4){\makebox(0,0)[lc]{\(3^4\)}}
\put(-5.2,3){\makebox(0,0)[rc]{\(27\)}}
\put(-4.8,3){\makebox(0,0)[lc]{\(3^3\)}}
\put(-5.2,2){\makebox(0,0)[rc]{\(9\)}}
\put(-4.8,2){\makebox(0,0)[lc]{\(3^2\)}}
\put(-5.2,1){\makebox(0,0)[rc]{\(3\)}}
\put(-5.2,0){\makebox(0,0)[rc]{\(1\)}}
\multiput(-5.1,0)(0,1){12}{\line(1,0){0.2}}
\put(-5,0){\line(0,1){11}}

\put(0,11.5){\makebox(0,0)[cb]{\(e=5\), \(m=8\)}}
\put(0,11){\circle*{0.2}}
\multiput(0,11)(-1,-1){2}{\line(1,-1){1}}
\multiput(0,11)(1,-1){2}{\line(-1,-1){1}}
\put(0,9){\circle*{0.2}}
\put(0,9){\line(0,-1){1}}
\put(0,8){\circle*{0.2}}
\multiput(0,8)(-1,-1){2}{\line(1,-1){1}}
\multiput(0,8)(1,-1){2}{\line(-1,-1){1}}
\put(0,6){\circle*{0.2}}
\multiput(0,6)(-1,-1){2}{\line(1,-1){1}}
\multiput(0,6)(1,-1){2}{\line(-1,-1){1}}
\put(0,4){\circle*{0.2}}
\multiput(0,4)(-1,-1){2}{\line(1,-1){1}}
\multiput(0,4)(1,-1){2}{\line(-1,-1){1}}
\put(0,2){\circle*{0.2}}

\put(1,7){\circle*{0.2}}
\multiput(-1,7)(-1,-1){2}{\line(1,-1){1}}
\multiput(-1,7)(1,-1){2}{\line(-1,-1){1}}
\multiput(1,7)(-1,-1){2}{\line(1,-1){1}}
\multiput(1,7)(1,-1){2}{\line(-1,-1){1}}

\multiput(-2,6)(-1,-1){2}{\line(1,-1){1}}
\multiput(-2,6)(1,-1){2}{\line(-1,-1){1}}

\multiput(-1,5)(-1,-1){2}{\line(1,-1){1}}
\multiput(-1,5)(1,-1){2}{\line(-1,-1){1}}
\multiput(1,5)(-1,-1){2}{\line(1,-1){1}}
\multiput(1,5)(1,-1){2}{\line(-1,-1){1}}
\multiput(3,5)(-1,-1){2}{\line(1,-1){1}}
\multiput(3,5)(1,-1){2}{\line(-1,-1){1}}

\multiput(4,4)(-1,-1){2}{\line(1,-1){1}}
\multiput(4,4)(1,-1){2}{\line(-1,-1){1}}

\multiput(1,3)(-1,-1){2}{\line(1,-1){1}}
\multiput(1,3)(1,-1){2}{\line(-1,-1){1}}
\put(1,1){\circle*{0.2}}
\multiput(3,3)(-1,-1){2}{\line(1,-1){1}}
\multiput(3,3)(1,-1){2}{\line(-1,-1){1}}

\put(2,6){\circle*{0.2}}
\multiput(2,6)(-1,-1){2}{\line(1,-1){1}}
\multiput(2,6)(1,-1){2}{\line(-1,-1){1}}
\put(2,4){\circle*{0.2}}
\multiput(2,4)(-1,-1){2}{\line(1,-1){1}}
\multiput(2,4)(1,-1){2}{\line(-1,-1){1}}
\put(2,2){\circle*{0.2}}
\multiput(2,2)(-1,-1){2}{\line(1,-1){1}}
\multiput(2,2)(1,-1){2}{\line(-1,-1){1}}
\put(2,0){\circle*{0.2}}

\put(9,10.5){\makebox(0,0)[cb]{\(e=5\), \(m=7\)}}
\put(9,10){\circle*{0.2}}
\multiput(9,10)(-1,-1){2}{\line(1,-1){1}}
\multiput(9,10)(1,-1){2}{\line(-1,-1){1}}
\put(9,8){\circle*{0.2}}
\put(9,8){\line(0,-1){1}}
\multiput(8,6)(-1,-1){2}{\line(1,-1){1}}
\multiput(8,6)(1,-1){2}{\line(-1,-1){1}}

\multiput(7,5)(-1,-1){2}{\line(1,-1){1}}
\multiput(7,5)(1,-1){2}{\line(-1,-1){1}}
\multiput(11,5)(-1,-1){2}{\line(1,-1){1}}
\multiput(11,5)(1,-1){2}{\line(-1,-1){1}}

\multiput(8,4)(-1,-1){2}{\line(1,-1){1}}
\multiput(8,4)(1,-1){2}{\line(-1,-1){1}}
\multiput(12,4)(-1,-1){2}{\line(1,-1){1}}
\multiput(12,4)(1,-1){2}{\line(-1,-1){1}}

\multiput(11,3)(-1,-1){2}{\line(1,-1){1}}
\multiput(11,3)(1,-1){2}{\line(-1,-1){1}}

\put(9,7){\circle*{0.2}}
\multiput(9,7)(-1,-1){2}{\line(1,-1){1}}
\multiput(9,7)(1,-1){2}{\line(-1,-1){1}}
\put(9,5){\circle*{0.2}}
\multiput(9,5)(-1,-1){2}{\line(1,-1){1}}
\multiput(9,5)(1,-1){2}{\line(-1,-1){1}}
\put(9,3){\circle*{0.2}}
\multiput(9,3)(-1,-1){2}{\line(1,-1){1}}
\multiput(9,3)(1,-1){2}{\line(-1,-1){1}}
\put(9,1){\circle*{0.2}}

\put(10,6){\circle*{0.2}}
\multiput(10,6)(-1,-1){2}{\line(1,-1){1}}
\multiput(10,6)(1,-1){2}{\line(-1,-1){1}}
\put(10,4){\circle*{0.2}}
\multiput(10,4)(-1,-1){2}{\line(1,-1){1}}
\multiput(10,4)(1,-1){2}{\line(-1,-1){1}}
\put(10,2){\circle*{0.2}}
\multiput(10,2)(-1,-1){2}{\line(1,-1){1}}
\multiput(10,2)(1,-1){2}{\line(-1,-1){1}}
\put(10,0){\circle*{0.2}}

\put(16,8.5){\makebox(0,0)[cb]{\(e=4\), \(m=6\)}}
\put(16,8){\circle*{0.2}}
\multiput(16,8)(-1,-1){2}{\line(1,-1){1}}
\multiput(16,8)(1,-1){2}{\line(-1,-1){1}}
\put(16,6){\circle*{0.2}}
\put(16,6){\line(0,-1){1}}
\multiput(15,4)(-1,-1){2}{\line(1,-1){1}}
\multiput(15,4)(1,-1){2}{\line(-1,-1){1}}

\put(16,5){\circle*{0.2}}
\multiput(16,5)(-1,-1){2}{\line(1,-1){1}}
\multiput(16,5)(1,-1){2}{\line(-1,-1){1}}
\put(16,3){\circle*{0.2}}
\multiput(16,3)(-1,-1){2}{\line(1,-1){1}}
\multiput(16,3)(1,-1){2}{\line(-1,-1){1}}
\put(16,1){\circle*{0.2}}

\put(17,4){\circle*{0.2}}
\multiput(17,4)(-1,-1){2}{\line(1,-1){1}}
\multiput(17,4)(1,-1){2}{\line(-1,-1){1}}
\put(17,2){\circle*{0.2}}
\multiput(17,2)(-1,-1){2}{\line(1,-1){1}}
\multiput(17,2)(1,-1){2}{\line(-1,-1){1}}
\put(17,0){\circle*{0.2}}

\multiput(18,3)(-1,-1){2}{\line(1,-1){1}}
\multiput(18,3)(1,-1){2}{\line(-1,-1){1}}

\end{picture}

\end{figure}

Figure
\ref{fig:BCFGroups}
shows several examples of normal lattices of \(3\)-groups \(G\)
with \textit{bicyclic and cyclic factors} of the central series.
They are located on coclass trees of coclass graphs \(\mathcal{G}(3,r)\)
\cite[p. 189 ff]{Ne}.

Here, the length of the rectangle of trailing diamonds is bigger than the width, \(m-1>e\),
the upper central series is different from the lower central series,
and the last lower central \(\gamma_{m-1}(G)\) is cyclic,
whence the parent \(\pi(G)=G/\gamma_{m-1}(G)\) is of the same coclass.
Such groups were called \textit{core groups} in
\cite{Ma}.
Concerning the principalization type \(\varkappa(K)\) of \(K\)
which coincides with the transfer kernel type (TKT) \(\varkappa(G)\) of \(G\), see
\cite{Ma1,Ma}.
Different TKTs can give rise to equal normal lattices.

\begin{example}
\label{exm:BCFGroups}
\(3\)-groups \(G\) of coclass \(3\le\mathrm{cc}(G)\le 4\).

\begin{itemize}
\item
Coclass \(\mathrm{cc}(G)=4\), class \(\mathrm{cl}(G)=7\):\\
a total of \(14\) complex quadratic fields, e. g.,\\
\(D=-159\,208\) with principalization type F.13,\\
\(D=-249\,371\) with principalization type F.12,\\
\(D=-469\,787\) with principalization type F.11,\\
\(D=-469\,816\) with principalization type F.7,\\
and a single real quadratic field of discriminant\\
\(D=8\,127\,208\) with principalization type F.13,\\
branch groups of depth \(1\),
visualized by Figure
\ref{fig:BCFGroups},
\(e=5\), \(m=8\).
\item
Coclass \(\mathrm{cc}(G)=4\), class \(\mathrm{cl}(G)=6\):\\
a single real quadratic field of discriminant\\
\(D=8\,491\,713\) with principalization type d\({}^\ast\).25,\\
mainline group,
visualized by Figure
\ref{fig:BCFGroups},
\(e=5\), \(m=7\).
\item
Coclass \(\mathrm{cc}(G)=3\), class \(\mathrm{cl}(G)=5\):\\
two real quadratic fields of discriminant\\
\(D=1\,535\,117\) with principalization type d.23,\\
\(D=2\,328\,721\) with principalization type d.19,\\
branch groups of depth \(1\),
visualized by Figure
\ref{fig:BCFGroups},
\(e=4\), \(m=6\).
\end{itemize}

\end{example}

\newpage


In Figure
\ref{fig:BFGroups}
we display numerous examples of normal lattices of \(p\)-groups \(G\)
with \textit{bicyclic factors} of the central series,
except the bottle neck \(\gamma_2(G)/\gamma_3(G)\).
They are located as vertices on the sporadic part \(\mathcal{G}_0(p,r)\)
of coclass graphs \(\mathcal{G}(p,r)\), outside of coclass trees,
\cite[Fig. 3.5, p. 439]{Ma}.

Here, the rectangle of trailing diamonds degenerates to a square with \(e=m-1\),
the upper central series is the reverse lower central series,
and thus the last lower central \(\gamma_{m-1}(G)\) is bicyclic,
whence the (generalized) parent \(\tilde\pi(G)=G/\gamma_{m-1}(G)\) is of lower coclass.
Such groups were called \textit{interface groups} in
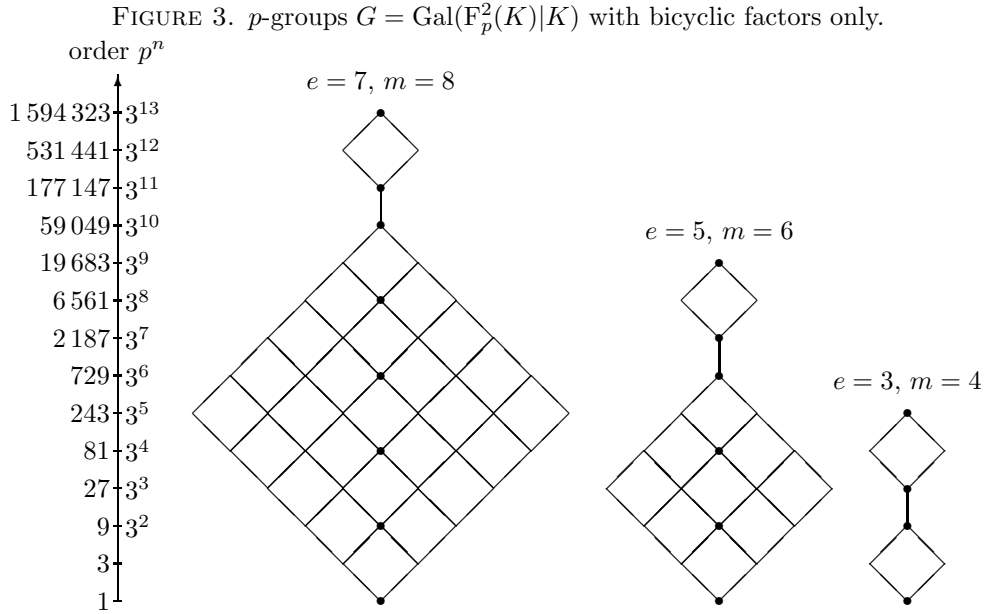
\cite{Ma}.

\begin{figure}[ht]
\caption{\(p\)-groups \(G=\mathrm{Gal}(\mathrm{F}_p^2(K)\vert K)\) with bicyclic factors only.}
\label{fig:BFGroups}


\setlength{\unitlength}{0.5cm}
\begin{picture}(24,15)(-7,0)

\put(-5,14.3){\makebox(0,0)[cb]{order \(p^n\)}}
\put(-5,13){\vector(0,1){1}}
\put(-5.2,13){\makebox(0,0)[rc]{\(1\,594\,323\)}}
\put(-4.8,13){\makebox(0,0)[lc]{\(3^{13}\)}}
\put(-5.2,12){\makebox(0,0)[rc]{\(531\,441\)}}
\put(-4.8,12){\makebox(0,0)[lc]{\(3^{12}\)}}
\put(-5.2,11){\makebox(0,0)[rc]{\(177\,147\)}}
\put(-4.8,11){\makebox(0,0)[lc]{\(3^{11}\)}}
\put(-5.2,10){\makebox(0,0)[rc]{\(59\,049\)}}
\put(-4.8,10){\makebox(0,0)[lc]{\(3^{10}\)}}
\put(-5.2,9){\makebox(0,0)[rc]{\(19\,683\)}}
\put(-4.8,9){\makebox(0,0)[lc]{\(3^9\)}}
\put(-5.2,8){\makebox(0,0)[rc]{\(6\,561\)}}
\put(-4.8,8){\makebox(0,0)[lc]{\(3^8\)}}
\put(-5.2,7){\makebox(0,0)[rc]{\(2\,187\)}}
\put(-4.8,7){\makebox(0,0)[lc]{\(3^7\)}}
\put(-5.2,6){\makebox(0,0)[rc]{\(729\)}}
\put(-4.8,6){\makebox(0,0)[lc]{\(3^6\)}}
\put(-5.2,5){\makebox(0,0)[rc]{\(243\)}}
\put(-4.8,5){\makebox(0,0)[lc]{\(3^5\)}}
\put(-5.2,4){\makebox(0,0)[rc]{\(81\)}}
\put(-4.8,4){\makebox(0,0)[lc]{\(3^4\)}}
\put(-5.2,3){\makebox(0,0)[rc]{\(27\)}}
\put(-4.8,3){\makebox(0,0)[lc]{\(3^3\)}}
\put(-5.2,2){\makebox(0,0)[rc]{\(9\)}}
\put(-4.8,2){\makebox(0,0)[lc]{\(3^2\)}}
\put(-5.2,1){\makebox(0,0)[rc]{\(3\)}}
\put(-5.2,0){\makebox(0,0)[rc]{\(1\)}}
\multiput(-5.1,0)(0,1){14}{\line(1,0){0.2}}
\put(-5,0){\line(0,1){13}}

\put(2,13.5){\makebox(0,0)[cb]{\(e=7\), \(m=8\)}}
\put(2,13){\circle*{0.2}}
\multiput(2,13)(-1,-1){2}{\line(1,-1){1}}
\multiput(2,13)(1,-1){2}{\line(-1,-1){1}}
\put(2,11){\circle*{0.2}}
\put(2,11){\line(0,-1){1}}
\put(2,10){\circle*{0.2}}
\multiput(2,10)(-1,-1){2}{\line(1,-1){1}}
\multiput(2,10)(1,-1){2}{\line(-1,-1){1}}
\put(2,8){\circle*{0.2}}
\multiput(2,8)(-1,-1){2}{\line(1,-1){1}}
\multiput(2,8)(1,-1){2}{\line(-1,-1){1}}
\put(2,6){\circle*{0.2}}
\multiput(2,6)(-1,-1){2}{\line(1,-1){1}}
\multiput(2,6)(1,-1){2}{\line(-1,-1){1}}
\put(2,4){\circle*{0.2}}
\multiput(2,4)(-1,-1){2}{\line(1,-1){1}}
\multiput(2,4)(1,-1){2}{\line(-1,-1){1}}
\put(2,2){\circle*{0.2}}
\multiput(2,2)(-1,-1){2}{\line(1,-1){1}}
\multiput(2,2)(1,-1){2}{\line(-1,-1){1}}
\put(2,0){\circle*{0.2}}

\multiput(1,9)(-1,-1){2}{\line(1,-1){1}}
\multiput(1,9)(1,-1){2}{\line(-1,-1){1}}
\multiput(3,9)(-1,-1){2}{\line(1,-1){1}}
\multiput(3,9)(1,-1){2}{\line(-1,-1){1}}

\multiput(0,8)(-1,-1){2}{\line(1,-1){1}}
\multiput(0,8)(1,-1){2}{\line(-1,-1){1}}
\multiput(4,8)(-1,-1){2}{\line(1,-1){1}}
\multiput(4,8)(1,-1){2}{\line(-1,-1){1}}

\multiput(-1,7)(-1,-1){2}{\line(1,-1){1}}
\multiput(-1,7)(1,-1){2}{\line(-1,-1){1}}
\multiput(1,7)(-1,-1){2}{\line(1,-1){1}}
\multiput(1,7)(1,-1){2}{\line(-1,-1){1}}
\multiput(3,7)(-1,-1){2}{\line(1,-1){1}}
\multiput(3,7)(1,-1){2}{\line(-1,-1){1}}
\multiput(5,7)(-1,-1){2}{\line(1,-1){1}}
\multiput(5,7)(1,-1){2}{\line(-1,-1){1}}

\multiput(-2,6)(-1,-1){2}{\line(1,-1){1}}
\multiput(-2,6)(1,-1){2}{\line(-1,-1){1}}
\multiput(0,6)(-1,-1){2}{\line(1,-1){1}}
\multiput(0,6)(1,-1){2}{\line(-1,-1){1}}
\multiput(4,6)(-1,-1){2}{\line(1,-1){1}}
\multiput(4,6)(1,-1){2}{\line(-1,-1){1}}
\multiput(6,6)(-1,-1){2}{\line(1,-1){1}}
\multiput(6,6)(1,-1){2}{\line(-1,-1){1}}

\multiput(-1,5)(-1,-1){2}{\line(1,-1){1}}
\multiput(-1,5)(1,-1){2}{\line(-1,-1){1}}
\multiput(1,5)(-1,-1){2}{\line(1,-1){1}}
\multiput(1,5)(1,-1){2}{\line(-1,-1){1}}
\multiput(3,5)(-1,-1){2}{\line(1,-1){1}}
\multiput(3,5)(1,-1){2}{\line(-1,-1){1}}
\multiput(5,5)(-1,-1){2}{\line(1,-1){1}}
\multiput(5,5)(1,-1){2}{\line(-1,-1){1}}

\multiput(0,4)(-1,-1){2}{\line(1,-1){1}}
\multiput(0,4)(1,-1){2}{\line(-1,-1){1}}
\multiput(4,4)(-1,-1){2}{\line(1,-1){1}}
\multiput(4,4)(1,-1){2}{\line(-1,-1){1}}

\multiput(1,3)(-1,-1){2}{\line(1,-1){1}}
\multiput(1,3)(1,-1){2}{\line(-1,-1){1}}
\multiput(3,3)(-1,-1){2}{\line(1,-1){1}}
\multiput(3,3)(1,-1){2}{\line(-1,-1){1}}

\put(11,9.5){\makebox(0,0)[cb]{\(e=5\), \(m=6\)}}
\put(11,9){\circle*{0.2}}
\multiput(11,9)(-1,-1){2}{\line(1,-1){1}}
\multiput(11,9)(1,-1){2}{\line(-1,-1){1}}
\put(11,7){\circle*{0.2}}
\put(11,7){\line(0,-1){1}}
\put(11,6){\circle*{0.2}}
\multiput(11,6)(-1,-1){2}{\line(1,-1){1}}
\multiput(11,6)(1,-1){2}{\line(-1,-1){1}}
\put(11,4){\circle*{0.2}}
\multiput(11,4)(-1,-1){2}{\line(1,-1){1}}
\multiput(11,4)(1,-1){2}{\line(-1,-1){1}}
\put(11,2){\circle*{0.2}}
\multiput(11,2)(-1,-1){2}{\line(1,-1){1}}
\multiput(11,2)(1,-1){2}{\line(-1,-1){1}}
\put(11,0){\circle*{0.2}}

\multiput(10,5)(-1,-1){2}{\line(1,-1){1}}
\multiput(10,5)(1,-1){2}{\line(-1,-1){1}}
\multiput(12,5)(-1,-1){2}{\line(1,-1){1}}
\multiput(12,5)(1,-1){2}{\line(-1,-1){1}}

\multiput(9,4)(-1,-1){2}{\line(1,-1){1}}
\multiput(9,4)(1,-1){2}{\line(-1,-1){1}}
\multiput(13,4)(-1,-1){2}{\line(1,-1){1}}
\multiput(13,4)(1,-1){2}{\line(-1,-1){1}}

\multiput(10,3)(-1,-1){2}{\line(1,-1){1}}
\multiput(10,3)(1,-1){2}{\line(-1,-1){1}}
\multiput(12,3)(-1,-1){2}{\line(1,-1){1}}
\multiput(12,3)(1,-1){2}{\line(-1,-1){1}}

\put(16,5.5){\makebox(0,0)[cb]{\(e=3\), \(m=4\)}}
\put(16,5){\circle*{0.2}}
\multiput(16,5)(-1,-1){2}{\line(1,-1){1}}
\multiput(16,5)(1,-1){2}{\line(-1,-1){1}}
\put(16,3){\circle*{0.2}}
\put(16,3){\line(0,-1){1}}
\put(16,2){\circle*{0.2}}
\multiput(16,2)(-1,-1){2}{\line(1,-1){1}}
\multiput(16,2)(1,-1){2}{\line(-1,-1){1}}
\put(16,0){\circle*{0.2}}

\end{picture}

\end{figure}

\begin{example}
\label{exm:BFGroups}
\(p\)-groups \(G\) with \(p\in\lbrace 3,5,7\rbrace\).

\begin{itemize}
\item
\(p=3\), coclass \(\mathrm{cc}(G)=6\), class \(\mathrm{cl}(G)=7\):\\
a single complex quadratic field of discriminant\\
\(D=-423\,640\) with principalization type F.12,\\
sporadic group,
visualized by Figure
\ref{fig:BFGroups},
\(e=7\), \(m=8\).
\item
\(p=3\), coclass \(\mathrm{cc}(G)=4\), class \(\mathrm{cl}(G)=5\):\\
a total of \(78\) complex quadratic fields, e. g.,\\
\(D=-27\,156\) with principalization type F.11,\\
\(D=-31\,908\) with principalization type F.12,\\
\(D=-67\,480\) with principalization type F.13,\\
\(D=-124\,363\) with principalization type F.7,\\
and a single real quadratic field of discriminant\\
\(D=8\,321\,505\) with principalization type F.13,\\
sporadic groups,
visualized by Figure
\ref{fig:BFGroups},
\(e=5\), \(m=6\).
\item
\(p=3\), coclass \(\mathrm{cc}(G)=2\), class \(\mathrm{cl}(G)=3\):\\
a total of \(936\) complex quadratic fields, e. g.,\\
\(D=-4\,027\) with principalization type D.10,\\
\(D=-12\,131\) with principalization type D.5,\\
and a total of \(140\) real quadratic fields, e. g.,\\
\(D=422\,573\) with principalization type D.10,\\
\(D=631\,769\) with principalization type D.5,\\
sporadic groups,
visualized by Figure
\ref{fig:BFGroups},
\(e=3\), \(m=4\).
\item
\(p=5\), coclass \(\mathrm{cc}(G)=2\), class \(\mathrm{cl}(G)=3\): see
\cite[Tbl. 3.13, p. 450]{Ma}.
\item
\(p=7\), coclass \(\mathrm{cc}(G)=2\), class \(\mathrm{cl}(G)=3\): see
\cite[Tbl. 3.14, p. 450]{Ma}.
\end{itemize}

\end{example}

\newpage


Figure
\ref{fig:SmallBCFGroups}
shows many examples of normal lattices of \lq\lq small\rq\rq\ \(p\)-groups \(G\)
with \textit{bicyclic and cyclic factors} of the central series.
They are located on coclass trees of coclass graphs \(\mathcal{G}(p,r)\)
\cite[Fig. 3.6--3.7, pp. 442--443]{Ma}.

\begin{figure}[ht]
\caption{Small \(p\)-groups \(G=\mathrm{Gal}(\mathrm{F}_p^2(K)\vert K)\) with bicyclic and cyclic factors.}
\label{fig:SmallBCFGroups}


\setlength{\unitlength}{0.5cm}
\begin{picture}(25,11)(-7,0)

\put(-5,10.3){\makebox(0,0)[cb]{order \(3^n\)}}
\put(-5,9){\vector(0,1){1}}
\put(-5.2,9){\makebox(0,0)[rc]{\(19\,683\)}}
\put(-4.8,9){\makebox(0,0)[lc]{\(3^9\)}}
\put(-5.2,8){\makebox(0,0)[rc]{\(6\,561\)}}
\put(-4.8,8){\makebox(0,0)[lc]{\(3^8\)}}
\put(-5.2,7){\makebox(0,0)[rc]{\(2\,187\)}}
\put(-4.8,7){\makebox(0,0)[lc]{\(3^7\)}}
\put(-5.2,6){\makebox(0,0)[rc]{\(729\)}}
\put(-4.8,6){\makebox(0,0)[lc]{\(3^6\)}}
\put(-5.2,5){\makebox(0,0)[rc]{\(243\)}}
\put(-4.8,5){\makebox(0,0)[lc]{\(3^5\)}}
\put(-5.2,4){\makebox(0,0)[rc]{\(81\)}}
\put(-4.8,4){\makebox(0,0)[lc]{\(3^4\)}}
\put(-5.2,3){\makebox(0,0)[rc]{\(27\)}}
\put(-4.8,3){\makebox(0,0)[lc]{\(3^3\)}}
\put(-5.2,2){\makebox(0,0)[rc]{\(9\)}}
\put(-4.8,2){\makebox(0,0)[lc]{\(3^2\)}}
\put(-5.2,1){\makebox(0,0)[rc]{\(3\)}}
\put(-5.2,0){\makebox(0,0)[rc]{\(1\)}}
\multiput(-5.1,0)(0,1){10}{\line(1,0){0.2}}
\put(-5,0){\line(0,1){9}}

\put(-2,9.5){\makebox(0,0)[cb]{\(e=3\), \(m=8\)}}
\put(-2,9){\circle*{0.2}}
\multiput(-2,9)(-1,-1){2}{\line(1,-1){1}}
\multiput(-2,9)(1,-1){2}{\line(-1,-1){1}}
\put(-2,7){\circle*{0.2}}
\put(-2,7){\line(0,-1){1}}
\put(-2,6){\circle*{0.2}}
\multiput(-2,6)(-1,-1){2}{\line(1,-1){1}}
\multiput(-2,6)(1,-1){2}{\line(-1,-1){1}}
\put(-2,4){\circle*{0.2}}

\put(-1,5){\circle*{0.2}}
\multiput(-1,5)(-1,-1){2}{\line(1,-1){1}}
\multiput(-1,5)(1,-1){2}{\line(-1,-1){1}}
\put(-1,3){\circle*{0.2}}

\put(0,4){\circle*{0.2}}
\multiput(0,4)(-1,-1){2}{\line(1,-1){1}}
\multiput(0,4)(1,-1){2}{\line(-1,-1){1}}
\put(0,2){\circle*{0.2}}

\put(1,3){\circle*{0.2}}
\multiput(1,3)(-1,-1){2}{\line(1,-1){1}}
\multiput(1,3)(1,-1){2}{\line(-1,-1){1}}
\put(1,1){\circle*{0.2}}

\put(2,2){\circle*{0.2}}
\multiput(2,2)(-1,-1){2}{\line(1,-1){1}}
\multiput(2,2)(1,-1){2}{\line(-1,-1){1}}
\put(2,0){\circle*{0.2}}

\put(5,8.5){\makebox(0,0)[cb]{\(e=3\), \(m=7\)}}
\put(5,8){\circle*{0.2}}
\multiput(5,8)(-1,-1){2}{\line(1,-1){1}}
\multiput(5,8)(1,-1){2}{\line(-1,-1){1}}
\put(5,6){\circle*{0.2}}
\put(5,6){\line(0,-1){1}}
\put(5,5){\circle*{0.2}}
\multiput(5,5)(-1,-1){2}{\line(1,-1){1}}
\multiput(5,5)(1,-1){2}{\line(-1,-1){1}}
\put(5,3){\circle*{0.2}}

\put(6,4){\circle*{0.2}}
\multiput(6,4)(-1,-1){2}{\line(1,-1){1}}
\multiput(6,4)(1,-1){2}{\line(-1,-1){1}}
\put(6,2){\circle*{0.2}}

\put(7,3){\circle*{0.2}}
\multiput(7,3)(-1,-1){2}{\line(1,-1){1}}
\multiput(7,3)(1,-1){2}{\line(-1,-1){1}}
\put(7,1){\circle*{0.2}}

\put(8,2){\circle*{0.2}}
\multiput(8,2)(-1,-1){2}{\line(1,-1){1}}
\multiput(8,2)(1,-1){2}{\line(-1,-1){1}}
\put(8,0){\circle*{0.2}}

\put(11,7.5){\makebox(0,0)[cb]{\(e=3\), \(m=6\)}}
\put(11,7){\circle*{0.2}}
\multiput(11,7)(-1,-1){2}{\line(1,-1){1}}
\multiput(11,7)(1,-1){2}{\line(-1,-1){1}}
\put(11,5){\circle*{0.2}}
\put(11,5){\line(0,-1){1}}
\put(11,4){\circle*{0.2}}
\multiput(11,4)(-1,-1){2}{\line(1,-1){1}}
\multiput(11,4)(1,-1){2}{\line(-1,-1){1}}
\put(11,2){\circle*{0.2}}

\put(12,3){\circle*{0.2}}
\multiput(12,3)(-1,-1){2}{\line(1,-1){1}}
\multiput(12,3)(1,-1){2}{\line(-1,-1){1}}
\put(12,1){\circle*{0.2}}

\put(13,2){\circle*{0.2}}
\multiput(13,2)(-1,-1){2}{\line(1,-1){1}}
\multiput(13,2)(1,-1){2}{\line(-1,-1){1}}
\put(13,0){\circle*{0.2}}

\put(16,6.5){\makebox(0,0)[cb]{\(e=3\), \(m=5\)}}
\put(16,6){\circle*{0.2}}
\multiput(16,6)(-1,-1){2}{\line(1,-1){1}}
\multiput(16,6)(1,-1){2}{\line(-1,-1){1}}
\put(16,4){\circle*{0.2}}
\put(16,4){\line(0,-1){1}}
\put(16,3){\circle*{0.2}}
\multiput(16,3)(-1,-1){2}{\line(1,-1){1}}
\multiput(16,3)(1,-1){2}{\line(-1,-1){1}}
\put(16,1){\circle*{0.2}}

\put(17,2){\circle*{0.2}}
\multiput(17,2)(-1,-1){2}{\line(1,-1){1}}
\multiput(17,2)(1,-1){2}{\line(-1,-1){1}}
\put(17,0){\circle*{0.2}}

\end{picture}

\end{figure}

\begin{example}
\label{exm:SmallBCFGroups}
Small \(p\)-groups \(G\) with \(p\in\lbrace 3,5,7\rbrace\).

\begin{itemize}
\item
\(p=3\), coclass \(\mathrm{cc}(G)=2\), class \(\mathrm{cl}(G)=7\):\\
a total of \(28\) complex quadratic fields, e. g.,\\
\(D=-262\,744\) with principalization type E.14,\\
\(D=-268\,040\) with principalization type E.6,\\
\(D=-297\,079\) with principalization type E.9,\\
\(D=-370\,740\) with principalization type E.8,\\
branch groups of depth \(1\),
visualized by Figure
\ref{fig:SmallBCFGroups},
\(e=3\), \(m=8\).
\item
\(p=3\), coclass \(\mathrm{cc}(G)=2\), class \(\mathrm{cl}(G)=6\):\\
two real quadratic fields, e. g.,\\
\(D=1\,001\,957\) with principalization type c.21,\\
mainline groups,
visualized by Figure
\ref{fig:SmallBCFGroups},
\(e=3\), \(m=7\).
\item
\(p=3\), coclass \(\mathrm{cc}(G)=2\), class \(\mathrm{cl}(G)=5\):\\
a total of \(383\) complex quadratic fields, e. g.,\\
\(D=-9\,748\) with principalization type E.9,\\
\(D=-15\,544\) with principalization type E.6,\\
\(D=-16\,627\) with principalization type E.14,\\
\(D=-34\,867\) with principalization type E.8,\\
and a total of \(21\) real quadratic fields, e. g.,\\
\(D=342\,664\) with principalization type E.9,\\
\(D=3\,918\,837\) with principalization type E.14,\\
\(D=5\,264\,069\) with principalization type E.6,\\
\(D=6\,098\,360\) with principalization type E.8,\\
branch groups of depth \(1\),
visualized by Figure
\ref{fig:SmallBCFGroups},
\(e=3\), \(m=6\).
\item
\(p=3\), coclass \(\mathrm{cc}(G)=2\), class \(\mathrm{cl}(G)=4\):\\
a total of \(54\) real quadratic fields, e. g.,\\
\(D=534\,824\) with principalization type c.18,\\
\(D=540\,365\) with principalization type c.21,\\
mainline groups,
visualized by Figure
\ref{fig:SmallBCFGroups},
\(e=3\), \(m=5\).
\item
\(p=5\), coclass \(\mathrm{cc}(G)=2\), class \(\mathrm{cl}(G)=5\): see
\cite[Tbl. 3.13, p. 450]{Ma}.
\item
\(p=7\), coclass \(\mathrm{cc}(G)=2\), class \(\mathrm{cl}(G)=5\): see
\cite[Tbl. 3.14, p. 450]{Ma}.
\end{itemize}

\end{example}

\newpage


\section{Final Remarks}
\label{s:FinalRemarks}

\begin{itemize}

\item
Among the \(2\,020\) complex quadratic fields
with \(3\)-class group of type \((3,3)\)
and discriminant in the range \(-10^6<D<0\),
the dominating part of \(1\,440\), that is \(71.29\,\%\),
has a second \(3\)-class group
with minimal defect of commutativity \(k=0\).
The remaining \(28.71\,\%\) have \(k=1\)
and TKTs G.16, G.19 and H.4.

\item
Among the \(2\,576\) real quadratic fields
with \(3\)-class group of type \((3,3)\)
and discriminant in the range \(0<D<10^7\),
a modest part of \(273\), i. e. \(10.6\,\%\),
has a second \(3\)-class group of coclass at least \(2\).
A dominating part of \(222\) among these \(273\) second \(3\)-class groups,
that is \(81.3\,\%\),
has minimal defect of commutativity \(k=0\),
whereas \(18.7\,\%\) have \(k=1\)
and TKTs b.10, G.16, G.19 and H.4.

\item
It should be pointed out that
the power-commutator presentations which we used for proving Theorem
\ref{thm:NormalLattice}
and its Corollaries are rudimentary,
since in fact they consist of commutator relations only.
Thus they define an isoclinism family of \(p\)-groups of fixed order,
rather than a single isomorphism class of \(p\)-groups.

On the other hand,
experience shows that the transfer kernel type (TKT) of a \(p\)-group
mainly depends on the power relations.
This explains why different TKTs frequently give rise to equal normal lattices.

\end{itemize}




\end{document}